\documentclass[12pt]{article}
\usepackage[margin=1in]{geometry}%

\usepackage[normalem]{ulem}
\usepackage{tikz}
\usetikzlibrary{arrows.meta,positioning}
\usetikzlibrary{arrows.meta,positioning}

\usetikzlibrary{positioning,fit,backgrounds}
\usepackage{hyperref}
\usepackage{amsmath,amsthm,enumitem,amssymb,mathtools}
\usepackage{graphicx}
\usepackage{float}
\usepackage{booktabs}
\usepackage{array}

\newcommand{\dist}{\mathrm{d}}

\newcommand{\SP}{\mathrm{SP}}
\newcommand{\length}{\mathrm{length}}
\newcommand{\HG}{\mathcal{H}}

\newcommand{\NB}{\mathrm{NB}}
\newcommand{\walk}{\mathrm{walk}}
\newcommand{\define}[1]{{\bf \boldmath{#1}}}

\usepackage{xcolor}

\theoremstyle{theorem}
\newtheorem{theorem}{Theorem}[section]
\newtheorem{corollary}[theorem]{Corollary}
\newtheorem{proposition}[theorem]{Proposition}
\newtheorem{lemma}[theorem]{Lemma}
\newtheorem{example}[theorem]{Example}
\newtheorem{remark}[theorem]{Remark}
\newtheorem{definition}[theorem]{Definition}

\title{Geometric characterisation of structural and regular equivalences in undirected (hyper)graphs}

\author{Marzieh Eidi$^{1,2}$ and Nina Otter$^{3,4}$\\
    $^1$ Center for Scalable Data Analytics and Artificial Intelligence (ScaDS.AI) and \\ 
    $^2$  Max Planck Institute for Mathematics in the Sciences (MPI MIS),\\ Leipzig, Germany \\
    $^3$ DataShape, Inria-Saclay, France \\
    $^4$ Laboratoire de Mathématiques d'Orsay, Université Paris-Saclay, France \\
}
\date{}
\begin{document}

\maketitle

\begin{abstract}
Similarity notions between vertices in a graph, such as structural and regular equivalence, are one of the main ingredients in clustering tools in complex network science. 
We generalise structural and regular equivalences for undirected hypergraphs and 
 provide a characterisation of structural and regular equivalences of undirected graphs and hypergraphs through neighbourhood graphs and Olivier-Ricci curvature. Our characterisation sheds new light on these similarity notions opening a new avenue for their exploration. 
 These characterisations also enable the construction of a possibly wide family of regular partitions, thereby offering a new route to a task that has so far been computationally challenging.
\end{abstract}
\tableofcontents

\section{Introduction}

Social Network Analysis  is a branch of social science characterised  by a structuralist approach that  puts the focus  on studying relationships between social actors rather than the actors themselves. 
The field has been developed since at least the 1930s, see \cite{freeman} for a historical  overview, and \cite{Wasserman_Faust_1994} for an overview of the main methods and questions. 
One of the main questions in the field is how to define a  position occupied by a given social actor in a  system. Different  equivalence relations  have  been introduced to capture what it means for two actors to occupy similar social positions. The first of these was the notion of \emph{structural equivalence}: 
two  actors are structurally equivalent  if they interact in the same way with the same  actors.   This notion is too rigid to be useful in practice, and thus one relaxes it through more flexible notions: on one hand one is interested in studying measures of approximate structural equivalence, capturing the extent to which two actors are structurally equivalent, and on the other hand one looks for more flexible equivalence relations, of which the most flexible one is given by regular equivalence.  Two actors are \emph{regularly equivalent}  if they interact with actors  who are themselves equivalent. 

Regular equivalence is a  more flexible notion than structural equivalence, however more difficult to compute in practice. The set of partitions satisfying the properties of a regular equivalence forms a lattice, and several algorithms for computing some elements of this lattice have been proposed, such as REGE \cite{rege1}, an algorithm that iteratively refines partitions and  returns  the maximal regular equivalence. More recently, algorithms that compute approximate regular partitions ahve been proposed in \cite{STT25}. See \cite{EB19} for an overview on different algorithms and challenges.
The situation becomes even more challenging in hypergraphs, which model  interactions between two or more actors. 
While recent years have seen advances in the development of generalisations of methods from graphs to hypergraphs \cite{bick,Zhou2007}, 
the development  of  structural and regular equivalence notions and algorithms for   hypergraphs are essentially absent.

In our work we propose a geometric characterisation of structural equivalence based on  Ollivier--Ricci curvature. Ollivier-Ricci curvature \cite{Ollivier2009} is a generalised notion of Ricci curvature for metric measure spaces, defined by comparing probability measures on neighbourhoods via the optimal transportation distance. Intuitively, curvature along the path connecting two vertices is positive if their neighbourhoods can be mapped to each other with little transportation cost. 
This makes curvature a natural tool to capture neighbourhood overlap in structural equivalence, and allows us to define different types of measures of  structural equivalence for hypergraphs, based on different types of random walks that capture different local patterns of  incidence. 

In a parallel line of investigation, we provide a way to compute some   regular partitions (i.e., partitions satisfying the property of being a regular equivalence) of graphs and hypergraphs based on their path-neighbourhood graphs, thus providing, to the best of our knowledge, the first algorithm to compute some  regular partitions for  hypergrahs. \\

Our main contributions are the following:

\begin{itemize}
\item We define notions of weak and strong structural and regular equivalence for hypergraphs, generalising structural and regular equivalence for graphs (see Definition \ref{D:str reg hyper}).
    \item We provide a curvature-based characterisation of structural equivalences in graphs (Theorem \ref{geom char eq}) and hypergraphs Theorem (\ref{T:hyper str}).  
    
    \item We propose a method to compute a class of regular partitions in graphs (Theorem \ref{T: char cc G_2} and Theorem \ref{T: reg G_n}) and hypergraphs (Corollary \ref{C:H_2 cc}(\ref{I: cc H_2})-(\ref{I: scc H_2}) and Corollary \ref{C: H_n reg}) through  connected components of the associated $n$-path neighbourhood graphs $G_n$. 
 
    \item In Section \ref{S:approx str eq} we develop a  curvature-based  measure to approximate structural equivalences for graphs and hypergraphs and relate it to cosine similarity, one of the main measures of approximate structural equivalence for graphs.
\end{itemize}

\subsection{Related work}
To the best of our knowledge, connections between curvature and structural or regular equivalences have not been studied before. 
A major part of our work consists in relating the connected component structure of neighborhood graphs based on paths with regular or structural classes. While neighborhood graphs based on paths have been studied since at least the  1980s, see \cite{simic1983graph, BGK19}, they are less well known than neighborhood graphs based on walks. The latter have been extensively studied in complex network science and adjacent fields because of their relationship with powers of the adjacency matrix of a graph. In particular, their connection with Ollivier Ricci curvature has been studied in \cite{Frank}.

Generalisations of regular equivalences for hypergraphs have only been considered, to the best of our knowledge, in the PhD thesis by  Borgatti \cite{borgatti89}, see also the review article \cite{EB93}, in which Borgatti developed an equivalent notion of regular equivalence, based on regular colorings.

To the best of our knowledge, algorithms to compute  regular partitions  based on neighbourhood graphs and curvature have not been considered before in the literature. Our approach is thus fundamentally different from existing approaches that compute iterative refinements of partitions \cite{rege1,STT25}.  

\section{Background on graphs}

Here we recall basic notions related to structural and regular equivalence of undirected graphs, approximate notions of these equivalences,  Ricci curvature on graphs and neighbourhood graphs. 

\begin{definition}
An \define{undirected graph} $G$ is a pair $(V,E)$ where $V$ is a non-empty set and $E\subset \mathcal{P}_{\leq 2}(V)$. We call the elements of $V$ \define{vertices} and the elements of $E$ \define{edges}. We call an edge $e\in E$ with cardinality $1$ a \define{loop}. Two vertices $x,y\in V$ are \define{neighbours} if there exists and edge $e$ such that $x,y\in e$. The \define{degree} $d_x$ of a vertex $x$ is the number of edges of which it is an element.

\end{definition}
We thus note that for us an undirected graph may contain loops,  it is  allowed to contain at most one edge between any two pairs of vertices. 

\begin{definition}
Let  $G=(V,E)$ be an undirected graph.
A \define{walk}  between two vertices $u,v\in V$ is a sequence $u=x_0 e_1 x_1e_2\dots x_{n-1}e_n x_n=v$ of  edges $e_1,\dots , e_n$ and vertices $x_0, x_1,\dots , x_{n-1}, x_n$ such that $x_i,x_{i+1}\in e_{i+1}$ for all $i=0,\dots , n-1$. A walk that starts and ends at the same vertex and in which no other vertices are repeated is called a \define{cycle}. We also call  a cycle containing $k$ distinct vertices a \define{$k$-cycle}. We then say that a cycle is \define{odd} (resp.\ \define{even}) if it is a $k$-cycle for $k$ odd (resp.\ even).
An undirected graph $G=(V,E)$ is \define{connected} if for any pair of vertices in $V$ there exists a walk between them. 
A \define{subgraph} of a graph $G=(V,E)$ is a graph $G'=(V',E')$ such that $V'\subset V$ and $E'\subset E$.
  A \define{connected component} of an undirected graph $G=(V,E)$ is a connected subgraph $G'$ of $G$ that is not  contained in any other connected subgraph of $G$. 
  A \define{complete} graph is a graph in which every two vertices are connected by an edge. 
\end{definition}

\subsection{Equivalence relations on graphs}\label{SS:EQ}
 \begin{definition}
 Given an undirected graph $G$, and an equivalence relation $\sim$ on its set of vertices, we say that the equivalence relation is a:

\begin{itemize}
\item \define{Structural equivalence} 
if for all vertices $a,b$ we have that  $a \sim b$ if and only if:\\
  $\{a,c\}$ is an edge in $G$ if and only if  $\{b,c\}$ is an edge in $G$.

\item \define{Regular equivalence}
if for all vertices $a,b$ we have that  $a\sim b$ if and only if: \\
if $\{a,c\}$ is an edge in $G$ then  there exists $d$ such that $\{b,d\}$ is an edge  and $c\sim d$.
\end{itemize}

Given a partition of the set of vertices of a graph $G$, we say that it is a \define{structural}, respectively \define{regular partition} partition if it satisfies the properties of a structural, resp.\ regular equivalence relation. Similarly, we call the blocks of such a partition \define{structural} resp.\ \define{regular classes} of $G$. 
\end{definition}

We note that there are at most two structural equivalences on a given graph $G$, of which one is the partition with only singleton blocks.
Every structural equivalence is a regular equivalence. On the other hand, not every regular equivalence is a structural equivalence, and in fact, there are in general many different partitions on the set of vertices of a given graph that satisfy the property of being a regular equivalence. 
The set of regular partitions of a given graph forms a lattice, ordered by inclusion. This lattice has a maximal element, which for an undirected graph is always given by the trivial partition with just one block, and a minimal element, given by the partition in which every block is a singleton \cite{WR83}.
Since maximal partitions play an important role in our work, we fix the following terminology:

\begin{definition}
We say that a partition of a set is \define{trivial} if it contains exactly one block.
\end{definition}

In practice, vertices in a graph are rarely structurally equivalent. 
Thus, one seeks to find measures that quantify the extent to which two vertices fail to be structurally equivalent. Many such notions have been introduced in the literature. Here we recall cosine similarity, which is one of the main notions used.

\begin{definition}\label{D:cosine sim}
Let $x,y$ be two vertices in an undirected graph $G$. 
The \define{cosine similarity} of $x$ and $y$ is defined as:

\[
\sigma_{x,y}=\frac{\eta_{xy}}{\sqrt{d_x}\sqrt{d_y}} \, ,
\]
where $n_{xy}$ is the number of common neighbours of $x$ and $y$.
\end{definition}

We note that while regular equivalences provide greater flexibility, notions of approximate regular equivalences are also studied. 
 We refer the reader to \cite[Section 12.4.4]{Wasserman_Faust_1994} for an overview of the main measures of regular equivalence.

\subsection{Ollivier-Ricci curvature}\label{SS:ORC}

In 2007 Ollivier \cite{Ollivier}  generalised Ricci curvature to metric spaces equipped with a random walk, of which graphs and hypergraphs are a particular example. Here we first provide the definition given by Ollivier, and we then provide the instantiation for graphs, while we will discuss the case of hypergraphs in Section \ref{S:hyper}.

\begin{definition}
Let $(X,d)$ be a complete and separable metric space togehter with a Borel $\sigma$-algebra.  A \define{random walk} $\mu$ on $X$ is a family of probability measures $\mu_x$ on $X$ for each $x \in X$, such that each measure $\mu_x$ depends measurably on $x$ and has finite first moment. Given two probability measures $\mu_x,\mu_y$, the \define{transportation distance} or \define{$1$-Wasserstein distance} between $\mu_x,\mu_y$ is 
  \[
W_1(\mu_x,\mu_y):=\inf_{\mathcal{E}\in \Pi(\mu_x,\mu_y)}\int_{(x,y)\in X\times X} d(x,y) d\mathcal{E}(x,y)\, 
  \]
 where  $\Pi(\mu_x,\mu_y)$ is the set of probability measures on $X \times X$ projecting to $\mu_x$ and $\mu_y$.
\end{definition}

\begin{definition} Let $(X,d)$ be a metric space equipped with a random walk $\mu:=(\mu_x)_{x\in X}$ and $x,y\in X$ be two distinct points.
The \define{Ollivier-Ricci curvature (ORC) of  $(x,y)$} is: 
   \[\kappa(x,y):=1-\frac{W_1(\mu_x,\mu_y)}{d(x,y)}\, .\] 
\end{definition}

To define Ollivier-Ricci curvature for undirected graphs, one usually takes as metric  the shortest-path  distance.

 \begin{definition}
Let $G=(V,E)$ be an undirected graph. 
A \define{path} between two vertices in $V$ is a walk in which the vertices and edges are all distinct. 
The \define{length} of a path $\gamma$ is the number $\mathrm{length}(\gamma)$ of edges it contains.   
The \define{shortest path} (or \define{combinatorial distance})  is the function $d_\SP\colon V\times V\to \mathbb{N}$ defined by
\[\dist_\SP(u,v)=\min_\gamma \length(\gamma) \, ,
\]
where the minimum is taken over all paths $\gamma$ between $u$ and $v$.
\end{definition}

We note that in the geometric study of graphs, less restrictive notions of paths, in which one may revisit certain vertices or edges, are often considered. While these notions are not suitable for the characterisation of structural or regular equivalences, since revisiting a vertex or edge would be adding redundant information to the equivalence relation, we recall their definitions here, and where appropriate, we will  explicitly discuss the differences that arise when adopting these alternative notions.

\begin{definition}
Let $G=(V,E)$ be an undirected graph. 
A \define{non-backtracking walk}  between two vertices $u,v\in V$ is a walk in which  we do not allow the pattern $x_i e_{i+1} x_{i+1} e_{i+1} x_i$, that is,   a walk cannot immediately return to a previously visited  vertex.
\end{definition}

To define Ollivier Ricci curvature of two arbitrary vertices in an unweighted graph, we define a random walk  by assigning a probability measure to each vertex based on randomly jumping from that vertex to its adjacent neighbours:

\begin{definition}[Ollivier Ricci Curvature on Graphs]
Let $G=(V,E)$ be a connected undirected graph with no loops.   
For each vertex $x\in V$, define the \define{neighbour measure} $\mu_x$ as

\[
\mu_x(z) =
\begin{cases}
\displaystyle \frac{1}{d_x} & \text{if $z$ is a neighbour of $x$}\, ,\\
0 & \text{otherwise.}
\end{cases}
\]

The \define{transportation distance} or \define{$1$-Wasserstein distance} between two neighbour   measures $\mu_x$ and $\mu_y$ is given by the solution to the following optimal transport problem:
\[
W_1(\mu_x, \mu_y)
    = \min_{\xi \in \Pi(\mu_x, \mu_y)}
      \sum_{u,v \in V} d(u,v)\, \xi(u,v),
\]
where $\Pi(\mu_x,\mu_y)$ is the set of all couplings of $\mu_x$ and $\mu_y$, that is,
\[
\Pi(\mu_x,\mu_y)
   = \bigl\{\, \xi:V\times V\to[0,1] \,\bigm|\,
      \sum_{v\in V} \xi(u,v)=\mu_x(u),\
      \sum_{u\in V} \xi(u,v)=\mu_y(v)
   \bigr\}.
\]

Then, for a pair of
vertices $x$ and $y$ in the same connected component, the \define{Ollivier--Ricci curvature} is defined by
\[
\kappa(x,y) := 1 - \frac{W_1(\mu_x, \mu_y)}{d_{SP}(x,y)} \, .
\]

\end{definition}

We note that the definition of Ollivier-Ricci curvature for graphs we give here has first been considered in \cite{Ollivier}, and \cite{LY2010}.
A ``lazy'' version of neighbour random walk, in which a random walk does not move from its starting vertex with a given probability $\alpha$ was  
introduced 
in \cite{lin11}. The version we consider here corresponds to the case $\alpha=0$. In \cite{lin11} the authors study a modification of the notion of Ricci curvature for graphs introduced in \cite{Ollivier,LY2010}, in which they consider the limit as $\alpha\to 1$.

It may help intuitively to think of  values of $\mu_x$ as \emph{sand} and those of $\mu_y$ as \emph{holes}. Since both $\mu_x$ and $\mu_y$ are probability measures, their total mass is equal to one. The central idea of optimal transport is then to determine the most efficient (i.e., lowest-cost) way to move the sand to the holes.

\subsection{Neighbourhood graphs}
A main ingredient in our  characterisations are neighbourhood graphs associated to undirected graphs or hypergraphs.  

\begin{definition}
Let $G=(V,E)$ be an undirected graph. Let $n\in\mathbb{N}$ be a positive integer. The \define{$n$th neighbourhood graph of $G$} is the graph $G_n=(V,E_n)$ with vertices those of $G$ and edges given by paths of length $n$ in $G$:
\[
E_n=\{\{u,v\}\mid \text{ there exist  a path between $u$ and $v$ of length $n$}  \}\, .
\]
\end{definition}

We note that while neighbourhood graphs based on walks might be a more familiar concept, since they are directly connected to powers of the adjacency matrix of a graph, neighbourhood graphs based on paths  have been studied since at least the 1980s, see  \cite{simic1983graph, BGK19}.

Since we will, where appropriate, discuss some differences in neighbourhood graphs based on different relaxations of paths, we give the definition here:

\begin{definition}
Let $G=(V,E)$ be an undirected graph. Let $n\in\mathbb{N}$ be a positive integer. The \define{$n$th non-bactracking neighbourhood graph of $G$} is the graph $G^\NB_n=(V,E^\NB_n)$ with vertices those of $G$ and edges given by non-backtracking walks of length $n$ in $G$:
\[
E^\NB_n=\{\{u,v\}\mid \text{ there exists  a non-backtracking walk between $u$ and $v$ of length $n$}  \}\, .
\]

The \define{$n$th walk neighbourhood graph of $G$} is the graph $G^\walk_n=(V,E^\walk_n)$ with vertices those of $G$ and edges given by walks of length $n$ in $G$:
\[
E^\walk_n=\{\{u,v\}\mid \text{ there exists  a walk between $u$ and $v$ of length $n$}  \}\, .
\]
\end{definition}

We note that neighbourhood graphs constructed from the three different notions of paths, namely,  walk, non-backtracking walk, and path generally do not coincide. By construction, the corresponding edge sets for any $n \ge 2$ satisfy the following inclusions:

\begin{equation}\label{eq:nbhdg}
E_n \subseteq E_n^\NB \subseteq E_n^\walk\, .
\end{equation}

\section{Geometric characterisation of structural and regular equivalences of undirected graphs}\label{S: geom char graphs}
In this section we study structural and regular equivalences of undirected graphs. We first start by discussing some properties of $n$-neighbourhood graphs which will be useful in the rest of the section.

\subsection{Connected components of $n$-neighbourhood graphs}

As $n$ increases, the  $n$-neighbourhood graph $G_n$ can become disconnected very quickly. 
For instance, once $n$ exceeds the diameter of the graph, this phenomenon is easily observed in the case of a star graph with $m$ edges: $G_2$ already has two connected components, while $G_3$ consists of 
$m+1$ isolated vertices. Here we first provide an upper bound on the number of connected components of $G_2$.

\begin{lemma}\label{lem:G2}
Let $G$ be a  connected graph. Then $G_2$ has at most two connected components. More precisely:
\begin{enumerate}
    \item If $G$ is non-bipartite (i.e.,  it contains an odd cycle), then $G_2$ is connected.
    \item If $G$ is bipartite with partition of the set of vertices $V=V_0\sqcup V_1$, then $G_2$ has exactly two connected components with vertex sets $V_0$ and $V_1$.
\end{enumerate}
\end{lemma}
\begin{proof}We note that if \(G\) is bipartite with partition of the vertex set given by \(V_0\) and \(V_1\), then by definition there are no edges between nodes belonging to \(V_0\) nor between nodes belonging to \(V_1\). Since \(G\) is connected, there exists a path between any two vertices in \(V_0\), and likewise between any two vertices in \(V_1\). Such paths necessarily alternate between the two partitions at each step and hence have even length. Conversely, paths between a given node in $V_0$ and a given node in $V_1$ must have odd length.   Consequently, \(V_0\) and \(V_1\) are the vertex sets of the two connected components of \(G_2\), and therefore \(G_2\) has exactly two components.

\medskip

If \(G\) is not bipartite, we show that \(G_2\) is connected, i.e., that any two vertices are connected by a path in \(G_2\). Let \(x\) and \(y\) be arbitrary vertices of \(G\). Since \(G\) is connected, there exists at least one path between them. If such a path has even length, then \(x\) and \(y\) lie in the same connected component of \(G_2\).

Suppose instead that no even-length path exists between \(x\) and \(y\) in \(G\). Then there must be a path of odd length between them; write this path as
\[
x = x_1, x_2, \ldots, x_{k}, x_{k+1} = y,
\]
where \(k\) is odd. Because there is no even-length path between \(x\) and \(y\), the given odd-length path cannot pass through any edge that is part of an odd cycle, since otherwise one could obtain an even-length path connecting $x$ and $y$, by traversing the remaining part of the cycle instead.
Since \(G\) is non-bipartite, it contains at least one odd cycle. 
Let  \(\alpha\) be a vertex on this cycle such that the path from \(x_k\) to \(\alpha\) passing through $y$ includes at least one edge of the cycle (see  illustration in Figure \ref{lem:G2} for an example).
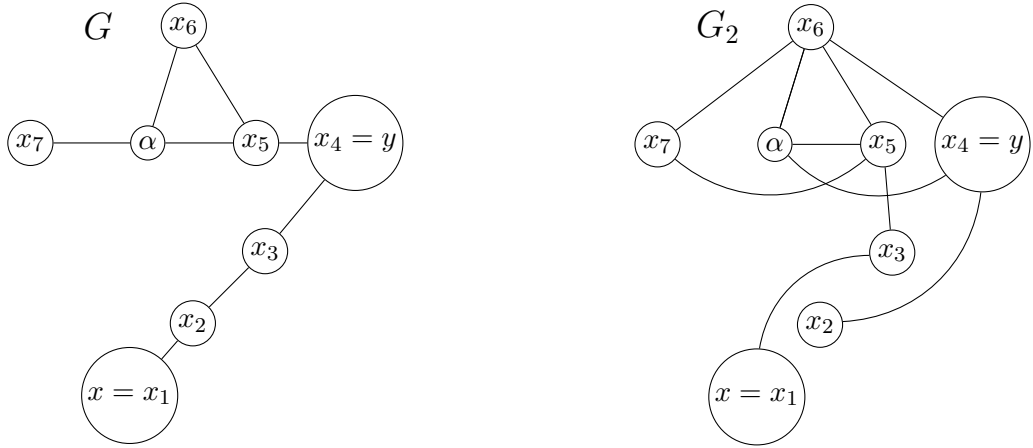
\begin{figure}[h]
\begin{center}

\begin{minipage}{0.54\textwidth}
\centering
\begin{tikzpicture}[
  scale=1.2,
  rotate=90,
  vertex/.style={circle,draw,fill=white,inner sep=1.6pt,minimum size=6pt,font=\small}
]

\node[vertex] (L1) at (-2.8, 1.3) {$x=x_1$};
\node[vertex] (L2) at (-2.0, 0.6) {$x_2$};
\node[vertex] (L3) at (-1.2,-0.2) {$x_3$};

\node[vertex] (B)  at ( 0.0,-1.2) {$x_4=y$};

\node[vertex] (V1) at ( 0.0,-0.1) {$x_5$};
\node[vertex] (V2) at ( 0.0, 1.1) {$\alpha$};
\node[vertex] (T)  at ( 0.0, 2.4) {$x_7$};

\node[vertex] (R)  at ( 1.3, 0.7) {$x_6$};

\draw (L1)--(L2)--(L3)--(B);
\draw (B)--(V1)--(V2)--(T);
\draw (V2)--(R)--(V1);

\node[left=14pt of R, font=\large] {$G$};

\end{tikzpicture}
\end{minipage}
\hfill
\begin{minipage}{0.45\textwidth}
\centering
\begin{tikzpicture}[
  scale=1.2,
  rotate=90,
  vertex/.style={circle,draw,fill=white,inner sep=1.6pt,minimum size=6pt,font=\small}
]

\node[vertex] (L1) at (-2.8, 1.3) {$x=x_1$};
\node[vertex] (L2) at (-2.0, 0.6) {$x_2$};
\node[vertex] (L3) at (-1.2,-0.2) {$x_3$};

\node[vertex] (B)  at ( 0.0,-1.2) {$x_4=y$};

\node[vertex] (V1) at ( 0.0,-0.1) {$x_5$};
\node[vertex] (V2) at ( 0.0, 1.1) {$\alpha$};
\node[vertex] (T)  at ( 0.0, 2.4) {$x_7$};

\node[vertex] (R)  at ( 1.3, 0.7) {$x_6$};

\draw (L1) to[bend left=40] (L3);
\draw (V1) to[bend left=40] (T);
\draw (L2) to[bend right=40] (B);
\draw (V1)--(V2)--(R);
\draw (R)--(T);
\draw (B)--(R);
\draw (V2) to[bend right=40] (B);
\draw (L3)--(V1);
\draw (V2)--(R)--(V1);

\node[left=14pt of R, font=\large] {$G_2$};

\end{tikzpicture}
\end{minipage}

\end{center}
\caption{An illustration of the construction in the argument of the proof of Lemma \ref{lem:G2} for $G$ non-bipartite. }
\end{figure}

Because the cycle has odd length, there exist two distinct paths between \(x_k\) and \(\alpha\), one of even length (call it path~1) and one of odd length (call it path~2). Hence, \(\alpha\) is connected to \(x_k\) in \(G_2\). Moreover, since \(x_k\) and \(y\) are adjacent in \(G\), \(\alpha\) is also connected to \(y\) in \(G_2\) (via a modified path \(2'\), obtained from path~2 by omitting the edge $\{x_k, y\}$ in $G$) and therefore, \(x_k\) and \(y\) lie in the same connected component of \(G_2\). Since \(x\) and \(x_k\) are connected by a path of even length, they are also in the same component. Hence, \(x\) and \(y\) belong to the same connected component of \(G_2\). As \(x\) and \(y\) were chosen arbitrarily, it follows that \(G_2\) is connected.

\end{proof}

One might wonder whether it is the case that the number of connected components of $G_n$ increases or decreases as $n$ varies. In general, neither is the case, as the next example shows. If $n$ passes the diameter of the graph, $G_n$ becomes disconnected and for large enough $n$, it is a collection of isolated vertices.

\begin{figure}
\[
\begin{tikzpicture}
\node at(0,-4) {$G=G_5$};
\node  at (0,0) {$\bullet$};
\node at (-1,-1) {$\bullet$};
\node at (1,-1) {$\bullet$};
\node at (-1,-2) {$\bullet$};
\node at (1,-2) {$\bullet$};
\node at (0,-3) {$\bullet$};
\path[-]
(0,0) edge (-1,-1)
(-1,-1) edge (-1,-2)
(-1,-2) edge (0,-3)
(0,0) edge (1,-1)
(1,-1) edge (1,-2)
(1,-2) edge (0,-3);
\end{tikzpicture}
\quad\quad
 \begin{tikzpicture}
 \node at(0,-4) {$G_2=G_4$};
\node (1) at (0,0) {$\bullet$};
\node (2) at (-1,-1) {$\bullet$};
\node at (1,-1) {$\bullet$};
\node at (-1,-2) {$\bullet$};
\node at (1,-2) {$\bullet$};
\node at (0,-3) {$\bullet$};
\path[-]
(0,0) edge (-1,-2)
(-1,-2) edge (1,-2)
(1,-2) edge (0,0)
(-1,-1) edge (1,-1)
(1,-1) edge (0,-3)
(0,-3) edge (-1,-1) ;
\end{tikzpicture}\quad\quad
\begin{tikzpicture}
\node at(0,-4) {$G_3$};
\node (1) at (0,0) {$\bullet$};
\node (2) at (-1,-1) {$\bullet$};
\node at (1,-1) {$\bullet$};
\node at (-1,-2) {$\bullet$};
\node at (1,-2) {$\bullet$};
\node at (0,-3) {$\bullet$};
\path[-]
(0,0) edge (0,-3)
(-1,-1) edge (1,-2)
(-1,-2) edge (1,-1);
\end{tikzpicture}\quad\quad
\begin{tikzpicture}
\node at(0,-4) {$G_6$};
\node (1) at (0,0) {$\bullet$};
\node (2) at (-1,-1) {$\bullet$};
\node at (1,-1) {$\bullet$};
\node at (-1,-2) {$\bullet$};
\node at (1,-2) {$\bullet$};
\node at (0,-3) {$\bullet$};
\end{tikzpicture}\\
\]
 \caption{Four different types of graphs that arise as $n$-neighbourhood graphs of $G$.
}\label{F:neighbourhood graphs}
\end{figure}
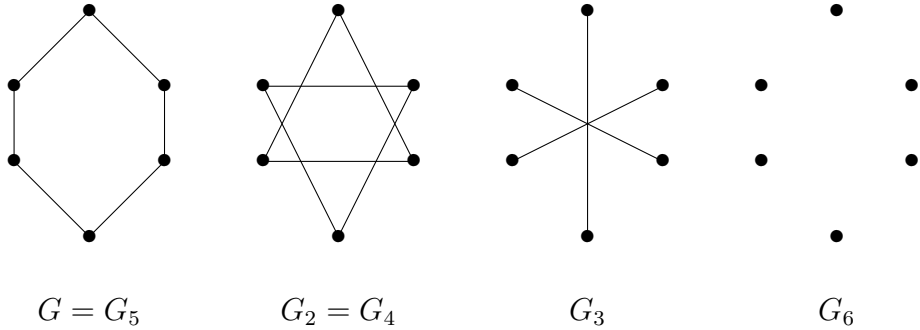

\begin{example}\label{E: connected components neighbourhood}
We consider the graph in Figure \ref{F:neighbourhood graphs} and we note that there are exactly four different types of graphs that arise as it $n$-neighbourhood graphs:

\[G_n=\begin{cases}
G,\; n=1 \text{ or } 5 \\
G_2,\; n=2 \text{ or } 4 \\
G_3,\; n=3 \\
G_6,\; n \ge 6  
\end{cases}
\]
Thus, we can see that the number of components of $G_n$ varies between  $1, 2, 3, 6$ as $n$ increases.
\end{example}

We note that, according to Equation \ref{eq:nbhdg}, the number of connected components satisfies
\[
\#\mathrm{Comp}\!\left(G_n^{\text{walk}}\right)
\;\le\;
\#\mathrm{Comp}\!\left(G_n^{\text{NB}}\right)
\;\le\;
\#\mathrm{Comp}\!\left(G_n\right).
\]

In other words, the $n$-neighbourhood graph defined by  walks is at least as connected as its non-backtracking counterpart, which in turn is at least as connected as the  $n$-neighbourhood graph.
We note that a key difference between the non-backtracking 
$n$-neighbourhood graph and the $n$-neighbourhood graph, lies in how they treat cycles. In the former case, a walk may traverse cycles of length at least three any number of times, whereas in the latter case, the path is not allowed to pass through a cycle more than once. Therefore for $n=2,3$ we have $G_n^{\text{NB}}=G_n$.

\subsection{Structural Equivalence}
Here we start to study the relationship between $n$-neighborhood graphs and structural partitions. For this we first fix some terminology:

\begin{definition}
    Let $n\geq 2$, and consider the $n$-neighbourhood graph $G_n$ of an undirected graph $G$. The classes \define{induced} by the connected components of $G_n$ are the blocks of a partition of $V$ given by vertices that belong to the same connected component. 
    \end{definition}

\begin{theorem} \label{geom char eq}
Let $G$ be a connected undirected graph without loops. 
 Then we have:

\begin{enumerate}[label=(\roman*)]
\item Two vertices $v,v'$ in $G$ are structurally equivalent if and only if the ORC of $(v,v')$ is equal to $1$. 
    \item Structural classes of $G$ are the vertex sets of  complete subgraphs of connected components of $G_2$. 
    \end{enumerate}
\end{theorem}
\begin{proof}
(i) By definition, two vertices x and y are structurally equivalent if and only if they have the same neighbours. This means that $\mu_x$ and $\mu_y$ defined above for $x$ and $y$ are the same and therefore the $1$-Wasserstein distance between them is equal to $0$. 
For the reverse, $\kappa(x,y)=1$ means that the $1$-Wasserstein distance between $\mu_x$ and $\mu_y$ is equal to $0$ and since $W_1$ is a metric, $\mu_x = \mu_y$ and thus $x$ and $y$ have exactly the same neighbours.

(ii) If $x$ and $y$ are structurally equivalent, they have exactly the same neighbours, and in particular they belong to the same connected component of $G_2$. We should show that for all  vertices in the same structural class, all  possible edges between them exist in $G_2$. Let us assume that  $x$ and $y$ are structurally equivalent.  If there is no edge between them, this means that there is no path of length $2$ between $x$ and $y$ in $G$, which contradicts the fact that they have at least one neighbour in common. 
\end{proof}

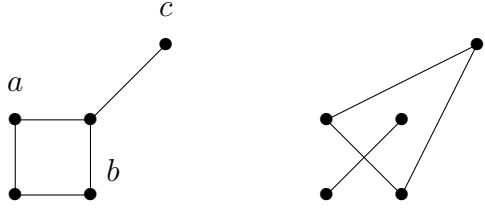
\begin{figure}
\begin{center}
\begin{tikzpicture}
\node at  (0,0) {$\bullet$};
\node at (1,0)[label={[shift={(0.3,-0.2)}]$b$}] {$\bullet$};
\node at (0,1) [label=$a$] {$\bullet$};
\node at (1,1) {$\bullet$};
\node at (2,2)[label=$c$] {$\bullet$};
\draw[-] (0,0)--(1,0)--(1,1)--(0,1)--(0,0);
\draw[-] (1,1)--(2,2);
\end{tikzpicture}\qquad\qquad
\begin{tikzpicture}
\node at (0,0) {$\bullet$};
\node at (1,0) {$\bullet$};
\node at (0,1) {$\bullet$};
\node at (1,1) {$\bullet$};
\node at (2,2) {$\bullet$};
\draw[-] (1,0)--(0,1)--(2,2)--(1,0);
\draw[-] (1,1)--(0,0);
\end{tikzpicture}
\end{center}
\caption{(Left) a graph $G$ and (right) its $2$-neighborhood graph $G_2$. While all connected components of $G_2$ are complete, the connected component on the nodes $a,b,c$  does not induce a structural class of $G$. }\label{F: G_2 complete non str}
\end{figure}
We note that the converse of statement (ii) is not true in general; namely, even if a connected component of $G_2$ is a complete graph, this doesn't necessarily imply that it induces a structural class,  as the example in Figure \ref{F: G_2 complete non str} shows.
However, by combining statements (i) and (ii) in Theorem  \ref{geom char eq} we obtain the following characterisation:

\begin{corollary}\label{C:compl G_2 str}
Complete subgraphs of connected components of $G_2$ induce  structural classes of $G$ if and only if the ORC  of every two pairs of vertices in the same subgraph is equal to $1$, and the subgraph is the maximal complete subgraph with respect to this property (i.e., there is no larger complete subgraph such that the ORC of all pairs of vertices in it is equal to $1$.)
\end{corollary} 

If we consider again the example from Figure \ref{F: G_2 complete non str}, we see that while the connected component on vertices $a,b,c$ does not induce a structural class, the complete subgraph given by the edge $\{a,b\}$ does. And in fact, this is the only pair of vertices from $\{a,b,c\}$
 having  ORC equal to $1$.

We note that based on the above characterisation, connected components of  $G_2$ provide limited information about structural classes of $G$. We will see next that if we instead consider the more general regular equivalences, connected components of $G_2$ always induce regular classes of $G$. 
Finally, let us observe that for two vertices $x$ and $y$ to be structurally equivalent, they cannot be neighbours and we need for them to have at least one neighbour in common,  thus they are necessarily in the same connected component of $G_2$. 
We will see in Section \ref{S:reg eq} how excluding this requirement enables us to go beyond $G_2$ and consider $G_n$ for general $n$.\\

\subsection{Regular Equivalence}\label{S:reg eq}
We now turn to regular equivalences. As mentioned in Section \ref{SS:EQ}, the set of regular partitions on a given vertex set $V$ of a graph $G$ forms a lattice. A question that has received a lot of interest is finding  efficient ways to compute elements of this lattice. In this section we explore the connection between elements of the lattice and connected components of  $n$-neighbourhood graphs. 
We note that in general non-trivial regular partitions of $G$ need not coincide with the connected components of $G_n$ for any $n$. For example, in a line graph with three edges, one may place the two endpoints in a single regular equivalence class and the two internal vertices in another one, even though this does not reflect the component structure of any $G_n$. Nevertheless, as shown in this section, there are several cases in which the component structure of $G_n$ can yield  a regular partition of $G$.
We start with the case of $G_2$.

\begin{theorem}\label{T: char cc G_2}
Connected components of $G_2$ induce regular classes of $G$. Moreover, if $G_2$ is not connected, i.e, has two connected components, these components are the blocks of a regular partition that is not properly contained in any  non-trivial partition. 
\end{theorem}
\begin{proof}
    For the first part, assume that $x$ and $y$ are two vertices in the same connected component of $G_2$.  Without loss of generality we may assume that $x$ and $y$ are connected by an edge in $G_2$. By definition of $G_2$, this means that in $G$ there is a vertex $\alpha$ such that $\alpha$ is simultaneously a neighbour of $x$ and $y$. Now, for an arbitrary neighbour $z$ of $x$ in $G$, if $z=\alpha$, then since $\alpha$ is also a neighbour of y,  the requirement for regular equivalence class is trivially satisfied. If  $z\neq \alpha$, then $z$ and $\alpha$ are connected by an edge in $G_2$, and thus they are in the same connected component. 
 This proves the first part of the  statement. 
    
We next observe that if $G_2$ is disconnected, then  the classes induced by the connected components of $G_2$ form a partition with two blocks. By construction, the only partition that is coarser than this is the trivial partition. 
\end{proof}

\begin{figure}[h]
\begin{center}
\begin{tikzpicture}
\node (1) at (0,0) {$\bullet$};
\node (2) at (1,0) {$\bullet$};
\node (3) at (2,0) {$\bullet$};
\node (4) at (3,0) {$\bullet$};
\node (5) at (4,0) {$\bullet$};
\node (6) at (5,0) {$\bullet$};
\node (7) at (3,-2) {$\bullet$};
\node (8) at (4,-2) {$\bullet$};
\node (9) at (5,-2) {$\bullet$};
\path[draw,-] (0,0)--(1,0)--(2,0)--(3,0)--(4,0)--(5,0)--(5,-2)--(4,-2)--(3,-2)--(3,0);
\end{tikzpicture}\qquad\qquad
\begin{tikzpicture}
\node (1) at (0,0) {\textcolor{violet}{$\bullet$}};
\node (2) at (1,0) {\textcolor{violet}{$\bullet$}};
\node (3) at (2,0) {\textcolor{orange}{$\bullet$}};
\node (4) at (3,0) {\textcolor{violet}{$\bullet$}};
\node (5) at (4,0) [label=$x$]{\textcolor{violet}{$\bullet$}};
\node (6) at (5,0)[label=$y$] {\textcolor{violet}{$\bullet$}};
\node (7) at (3,-2) {\textcolor{violet}{$\bullet$}};
\node (8) at (4,-2) {\textcolor{violet}{$\bullet$}};
\node (9) at (5,-2)  {\textcolor{orange}{$\bullet$}};
\path[draw,-] (0,0)--(1,0)--(2,0)--(3,0)--(4,0)--(5,0)--(5,-2)--(4,-2)--(3,-2)--(3,0);
\end{tikzpicture}
\end{center}
\caption{A graph $G$ (left), and  the partition of its set of nodes into two  blocks induced by the two connected components of  $G_4$, where we color vertices in the same block by the same color (right). The two connected components of $G_4$ do not induce regular classes on $G$: the nodes $x$ and $y$ belong to the same class, however, while $y$ is connected through an edge to an orange node, there is no edge from $x$ to  an orange node.}\label{F:G_4 not reg cla}
\end{figure}
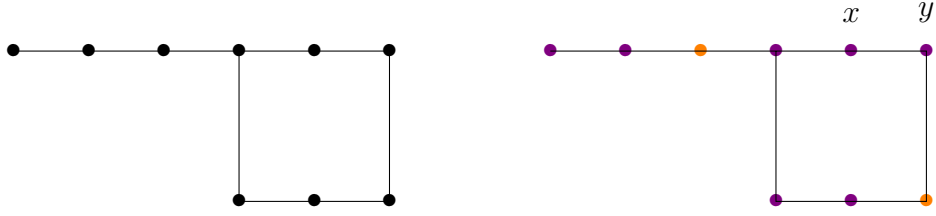

We note that the previous result does not directly generalise to $G_n$ for arbitrary $n$, namely, it is not true in general that connected components of $G_n$ give regular classes, not even in the case that $G_n$ has exactly two connected components, as the  example in Figure \ref{F:G_4 not reg cla} shows. 

Nonetheless, we can generalise the result under the additional assumption that no vertex has degree $1$.
\begin{theorem}\label{T: reg G_n} Let $G$ be an undirected graph. 
    If all vertices in $G$ have degree at least $2$ then the connected components of $G_n$ induce regular classes of $G$. 
\end{theorem}
\begin{proof}
    Suppose nodes $a$ and $b$ are in the same connected component of $G_n$ and without loss of generality, assume that they are connected by an edge in $G_n$. Therefore, there is a path of length $n$ between them in $G$ such as: $a, x_1,\dots , x_{n-1}, b$. Now for an arbitrary neighbour of $a$ such as $\alpha$, if $\alpha \neq x_1$, then $\{\alpha, x_{n-1}\}$ is an edge in $G_n$ where $x_{n-1}$ is a neighbour of $b$ in $G$. If 
    $\alpha = x_1$, since $d_b \geq 2$ there exists a neighbhour $\beta $ of $b$ such that $\beta \neq x_{n-1}$ and $\{x_1, \beta\}$ is an edge in $G_n$. By definition of regular equivalence, this proves the claim.    
\end{proof}

\begin{figure}[h]
\[
\begin{tikzpicture}
\node  at (0,-5) {$G$};
\node [label=a] at (0,0) {$\bullet$};
\node at (0,-1) {$\bullet$};
\node at (1,-2) {$\bullet$};
\node at (-1,-2) {$\bullet$};
\node at (1,-3) {$\bullet$};
\node at (-1,-3) {$\bullet$};
\node at (0,-4) {$\bullet$};   
\path[-] (0,0) edge (0,-1)
(0,-1) edge (-1,-2)
(-1,-2) edge (-1,-3)
(-1,-3) edge (0,-4)
(0,-4) edge (1,-3)
(0,-1) edge (1,-2)
(1,-2) edge (1,-3);
\end{tikzpicture}\quad \quad\quad
\begin{tikzpicture}
\node at (0,-5) {$G_3$};
\node[label=a] at (0,0) {$\bullet$};
\node at (0,-1) {$\bullet$};
\node at (1,-2) {$\bullet$};
\node at (-1,-2) {$\bullet$};
\node at (1,-3) {$\bullet$};
\node at (-1,-3) {$\bullet$};
\node at (0,-4) {$\bullet$};   
\path[-] (0,0) edge (-1,-3)
(0,0) edge (1,-3)
(1,-2) edge (-1,-3)
(-1,-2) edge (1,-3)
(0,-1) edge (0,-4)
;
\end{tikzpicture}
\]
\caption{An example of a graph $G$ with vertices of degree $1$ and in which the components of $G_3$ give regular partitions.}\label{E:G_3 comp}
\end{figure}
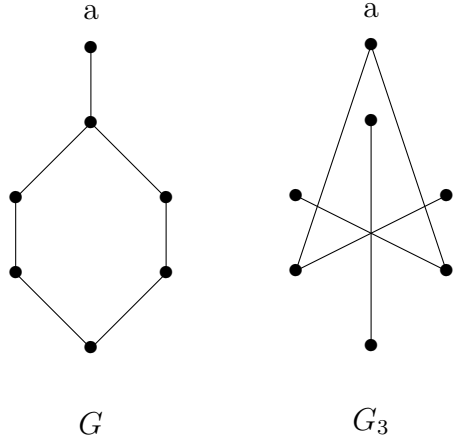

\begin{remark}
    The condition on the degree of vertices is not a necessary condition for the connected components of $G_n$ to induce regular classes of $G$.
For an example, consider the graph in Figure \ref{E:G_3 comp}. The components of $G_3$ induce a regular partition on $G$, however, vertex $a$ has degree $1$ in $G$.
\end{remark}


Now, focusing on the second neighbourhood graph $G_2$, we ask the following question:
\begin{equation}(*)\notag\label{Q}\text{
 Can we get other non-trivial regular classes from the connected components of $G_2$?}
 \end{equation}
 One answer can be obtained by relaxing the maximality condition in Corollary \ref{C:compl G_2 str}.

\begin{corollary}\label{T:subc}
    Complete subgraphs of the connected components of $G_2$ such that the curvature of all pairs of  vertices in each subgraph is equal to $1$ in $G$ induce regular classes of $G$. 
\end{corollary}
\begin{proof}
     Suppose that $x$ and $y$ are two vertices in a complete  subgraph of a connected component of $G_2$ such that curvature of $(x,y)$ in $G$ is equal to $1$. This means that their neighbours in $G$ coincide, and therefore the subgraph induces a regular equivalence class.
    
\end{proof}

\begin{figure}[h]
\centering

\usetikzlibrary{arrows.meta,calc,positioning}
\tikzset{
  vertex/.style = {circle, draw, thick, inner sep=1pt, minimum size=6pt, fill=white},
  vblue/.style  = {vertex, fill=blue!60},
  vgreen/.style = {vertex, fill=green!60!black},
  lab/.style    = {font=\small},
  edge/.style   = {line width=1pt},
  faded/.style  = {line width=0.9pt, opacity=0.95},
  arrow/.style  = {-{Latex[length=2.6mm,width=2mm]}, line width=0.9pt}
}

\begin{tikzpicture}[scale=0.8, transform shape]

\newcommand{\DrawG}[5]{%
  \coordinate (a) at (-2,0);
  \coordinate (e) at ( 2,0);
  \coordinate (b) at ( 0,1.5);
  \coordinate (d) at ( 0,-1.5);
  \coordinate (c) at ( 0,0);

  \draw[faded] (a)--(b)--(e)--(d)--(a);
  \draw[faded] (a)--(c)--(e);

  \node[#1] (A) at (a) {};
  \node[#2] (B) at (b) {};
  \node[#3] (C) at (c) {};
  \node[#4] (D) at (d) {};
  \node[#5] (E) at (e) {};

  \node[lab, left=1pt of A] {$a$};
  \node[lab, above=1pt of B] {$b$};
  \node[lab, above=1pt of C] {$c$};
  \node[lab, below=1pt of D] {$d$};
  \node[lab, right=1pt of E] {$e$};
}


\begin{scope}[shift={(-7,1.2)}]
  \node[lab] at (0,2.4) {$G$}; 
  \DrawG{vertex}{vertex}{vertex}{vertex}{vertex}
\end{scope}

\draw[arrow] (-4,1.2) -- (-3,1.2);

\begin{scope}[shift={(0,1.2)}]
  \node[lab] at (-3.2,1.2) {$G_2$}; 

  \coordinate (b2) at (0,1.2);
  \coordinate (c2) at (1.6,0);
  \coordinate (d2) at (-1.6,0);

  \draw[edge] (b2) to (c2);
  \draw[edge] (c2) to (d2);
  \draw[edge] (d2) to (b2);

  \node[vertex] at (b2) {};
  \node[vertex] at (c2) {};
  \node[vertex] at (d2) {};
  \node[lab, above=2pt] at (b2) {$c$};
  \node[lab, right=2pt] at (c2) {$d$};
  \node[lab, left=2pt]  at (d2) {$b$};

  \coordinate (a2) at (3,0.2);
  \coordinate (e2) at (5,0.2);
  \draw[edge] (a2)--(e2);
  \node[vertex] at (a2) {};
  \node[vertex] at (e2) {};
  \node[lab, below=1pt] at (a2) {$a$};
  \node[lab, below=1pt] at (e2) {$e$};
\end{scope}


\begin{scope}[shift={(-9,-3.5)}]
  \DrawG{vgreen}{vblue}{vblue}{vblue}{vgreen}
\end{scope}

\begin{scope}[shift={(-3.5,-3.5)}]
  \DrawG{vgreen}{vblue}{vblue}{vertex}{vgreen}
\end{scope}

\begin{scope}[shift={(2,-3.5)}]
  \DrawG{vgreen}{vertex}{vblue}{vblue}{vgreen}
\end{scope}

\begin{scope}[shift={(7.5,-3.5)}]
  \DrawG{vgreen}{vblue}{vertex}{vblue}{vgreen}
\end{scope}

\end{tikzpicture}

\caption{Top: a graph \(G\) (left) and its \(G_2\) (right), which splits into a triangle on \(\{b,c,d\}\) and an edge \(\{a,e\}\) (right). Bottom: Complete subgraphs of the  components of $G_2$, such that the curvature of every pair of vertices is equal to $1$ (in $G$), induce regular classes of $G$; vertices with the same color are in the same regular class. }\label{F:G_2 subc}
\end{figure}
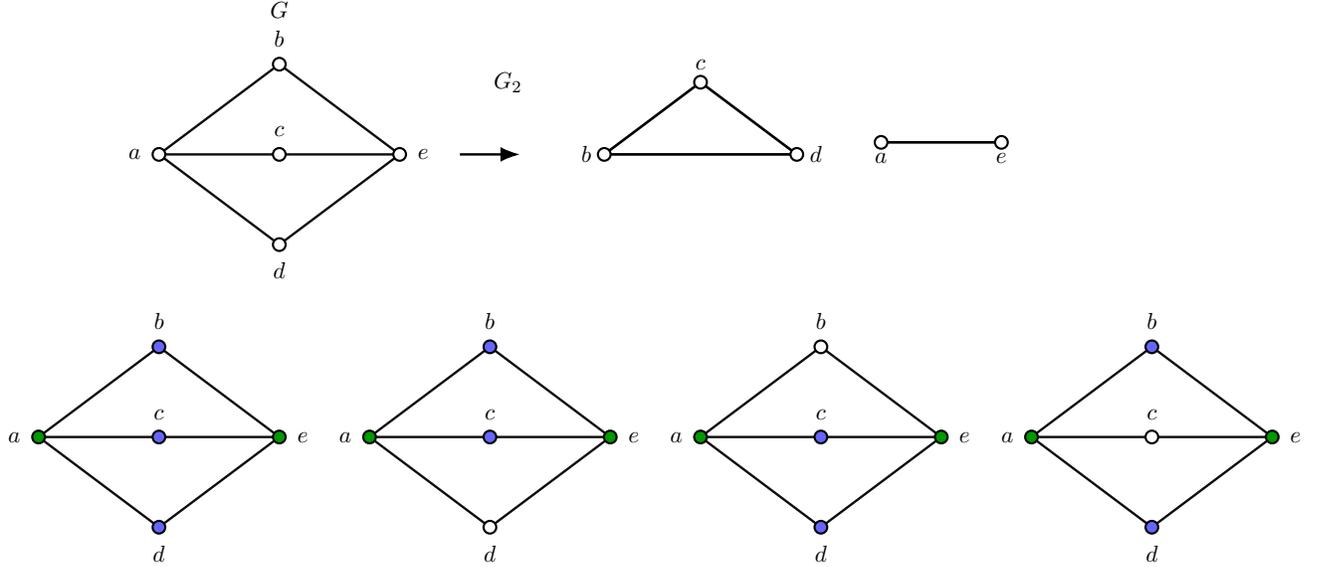

Thus, by Corollary \ref{T:subc} we can obtain a regular partition of $G$ by taking blocks induced by complete subgraphs of the connected components of $G_2$, while the remaining vertices are each assigned to a singleton block. We provide an example  in Figure \ref{F:G_2 subc}.
We note that, in the resulting partition, no two vertices belonging to the same class are connected by an edge in 
$G$. In fact, such partitions can be seen  as being  refinements of a given structural partition (in Figure \ref{F:G_2 subc} this would be the first partition on the bottom left).
 
The following proposition presents another answer to our question $(*)$ by a slight modification of  Theorem \ref{T: char cc G_2}.  
\begin{proposition}\label{T:G_2'}
Let $G'$ be the graph obtained from G by removing at most one edge  from each triangle in G and let $G_2'$ be the $2$-neighbourhood graph of $G'$. Then the connected components of $G_2'$ induce regular classes for the original graph $G$.
\end{proposition}
\begin{proof}
    We should show that if $\{x,y\}$ is an edge in $G_2'$ and $z$ is an arbitrary neighbour of $x$ in $G$, a neighbour of $y$ (such as $z'$) in $G$ exists such that $z$ and $ z'$  are in the same connected component of $G_2'$.
    By construction,  if $\{x,y\}$ is an edge in  $G_2'$ , there is $\alpha$ such that $\{x,\alpha\}$ and $\{\alpha,y\}$ are edges in $G'$ (and $G$). Now there are two possibilities: either there is an edge between $x$ and $y$ in $G$ (namely $\alpha, x, y$ constitute a triangle) or there is no such edge (there is no such triangle):
    \begin{itemize}
        \item If there is an edge between $x$ and $y$ in $G$, an arbitrary neighbour of $x$ could be $\alpha$, $y$ or neither of the two. In the first and second case, the condition of a regular class is trivially satisfied. In the last case, denote by $\alpha'$ the neighbour of $x$. Then there is an edge between $\alpha$ and $\alpha'$ in $G_2'$ and thus $\alpha$ and $\alpha'$ are int he same connected component of $G_2'$. Since $\alpha$ is a neighbour of $y$ in $G$, this prove the claim.
        \item  If there is no edge between $x$ and $y$ in $G$, then an arbitrary neighbour of $x$ could be $\alpha$ or a vertex $\alpha'$ different from $\alpha$ and $y$. The proof is similar to the one given in previous item, and we omit it here. 
      
    \end{itemize}
      
\end{proof}

We have seen in Theorem \ref{T: char cc G_2} that the connected components of $G_2$ induce regular partitions of $G$. Now, if $G$ is non-bipartite, by Lemma \ref{lem:G2} we know that this regular partition is the trivial partition.  Proposition \ref{T:G_2'} then provides a way to obtain non-trivial regular partitions from $G_2$ in the case that the odd cycles are given by triangles. 
We can thus ask whether it is possible to obtain a similar  statement for arbitrary odd cycles. It is indeed possible, however at the cost of adding some further restrictions on the edges that one may delete.
\begin{proposition}
Let $G$ be a connected undirected graph with no vertices of degree $1$.
Let $G'$ be the graph obtained from $G$ by removing at most one edge  from $k$-cycles ($k \ge 3$) in $G$ 
and such that no deletion of edges may result in vertices of degree $1$ in $G'$ (in other words, there are no vertices of degree $1$ in $G'$).
Let $G_{k-1}'$ be the $(k-1)$-neighborhood graph of $G'$. Then the connected components of $G_{k-1}'$ induce regular classes for the original graph $G$.
\end{proposition}

\begin{proof}
Let $\{x,y\}$ be an edge in $G'_{k-1}$. Let $\alpha$ be a neighbour of $x$ in $G$. We need to show that there exists a neighbour $\beta$ of $y$ in $G$ such that $\alpha$ and $\beta$ are in the same connected component of $G'_{k-1}$. By construction, there is a path of length $k-1$ between $x$
 and $y$, such as $x,x_1, \dots, x_{k-2}, y$.  There are two cases: (i) If $\alpha=x_1$, then since $y$ has degree at least $2$ in $G'$, there exists a neighbour $\beta$ of $y$ different from $x_{k-2}$ (in $G'$) and thus $\alpha,\beta$ are in the same connected component of $G_{k-1}'$. (ii) If $\alpha\ne x_1$ then $\alpha$ and $x_{k-2}$ are in the same connected component of $G_{k-1}'$.
 \end{proof}
\section{The case of undirected hypergraphs}\label{S:hyper}

To the best of our knowledge, the notions of structural and regular equivalence   have not yet been generalised to higher-order combinatorial structures such as hypergraphs. 
While both notions can be generalised to hypergraphs in a straightforward manner (see Definition \ref{D:str reg hyper}), the main difficulty  resides in how to compute partitions that satisfy the conditions of regularity. 
A framework relying on the theory of coalgebras that provides efficient computation of regular equivalences for a particular class of directed hypergraphs is being developed in current work by one of the authors \cite{MOR25}. Here we extend some of our geometric characterisations for undirected hypergraphs.

\begin{definition}
An \define{undirected hypergraph} $\HG$ is a 
pair $(V,H)$ where $V$ is a non-empty set and $H\subset \mathcal{P}(V)$ is a collection of subsets of $V$. We call the elements of $V$ \define{vertices of $\HG$} and the elements of $H$ \define{hyperedges of $\HG$}. The \define{cardinality} $|e|$ of a hyperedge $e$ is the number of elements of $e$. We call a hyperedge with cardinality equal to $2$ an \define{edge} and one with cardinality equal to $1$ a \define{loop}.  Two vertices $u,v$ are \define{neighbours} if there exists a hyperedge $e$ such that $u,v\in e$.  
The \define{degree} $d_v$ of a vertex $v\in V$ is the number of hyperedges to which it belongs.

\end{definition}

We next introduce two different notions of structural and regular equivalence, In the weak version we only consider the neighbours of two given vertices, disregarding the types of ties that connect them to the neighbouring vertices. In the strong version, we additionally take into account the types of ties.

\begin{definition}\label{D:str reg hyper}
Consider an undirected hypergraph $\HG=(V,H)$ together with an equivalence relation $\sim$ on its set of vertices. We say that $\sim$ is  a:
\begin{itemize}
\item \define{weak structural equivalence}  if for all vertices $u,v\in V$ we have that $u\sim v$ if and only if: $\{u,a\}\subset e$ for some hyperedge $e$ exactly when   $\{v,a\}\subset e'$ for some hyperedge $e'$. In other words, $u,v$ have exactly the same neighbours.
\item \define{strong  structural equivalence}  if for all vertices $u,v\in V$ we have that $u\sim v$ if and only if: $\{u,a\}\subset e$ for some hyperedge $e$ exactly when   $\{v,a\}\subset e'$ for some hyperedge $e'$ (i.e., $\sim$ is a weak structural equivalence) and additionally we have that $|e|=|e'|$. In other words, $u,v$ are connected to the same  neighbours through hyperedges of the same cardinality.
\item \define{weak regular equivalence} if for all vertices $u,v$ we have $u\sim v$ if and only if: if $\{u,a\}\subset e$ for a hyperedge $e$ then there exists $a'\in V$ and $e'\in H$ such that $\{v,a'\}\subset e'$ and $a\sim a'$. In other words, $u,v$ have neighbours who are themselves equivalent. 
\item \define{strong regular equivalence} 
if for all vertices $u,v$ we have $u\sim v$ if and only if: if $\{u,a\}\subset e$ for a hyperedge $e$ then there exists $a'\in V$ and $e'\in H$ such that $\{v,a'\}\subset e'$ and $a\sim a'$ (i.e., $\sim$ is a weak regular equivalence) and additionally $|e|=|e'|$. In other words, $u,v$ are connected to neighbours who are themselves equivalent through hyperedges of the same cardinality. 
\end{itemize}
We say that two vertices $u,v\in V$ are \define{strongly} (resp.\ \define{weakly}) \define{structurally} (resp.\ \define{regularly}) \define{equivalent} if there exists an equivalence relation on $V$ that is a strong (resp.\ weak) structural (resp.\ regular) equivalence. 

\end{definition}

\begin{remark}
We note that if two vertices are strongly structurally equivalent, then they necessarily have the same degree. Conversely, if two vertices have the same neighbours and same degree, it is not necessarily the case that they are connected to the same neighbours through hyperedges of the same cardinality. One could thus envision this to give us an intermediate notion of structural equivalence. For  a more in depth discussion of several different possible  notions of structural equivalences for hypergraphs see Remark \ref{R:hierarchy}.
\end{remark}

Now a main question arises: can we generalize the  results from Section \ref{S: geom char graphs} to hypergraphs? To answer this question, we need to first generalise some notions.

\begin{definition}
Let  $\HG=(V,E)$ be an undirected hypergraph.
A \define{walk}  between two vertices $u,v\in V$ is a sequence $u=x_0 e_1 x_1e_2\dots x_{n-1}e_n x_n=v$ of  hyperedges $e_1,\dots , e_n$ and vertices $x_0, x_1,\dots , x_{n-1}, x_n$ such that $x_i,x_{i+1}\in e_{i+1}$ for all $i=0,\dots , n-1$. A walk that starts and ends at the same vertex is called a \define{cycle}. A \define{path} is a walk in which all the vertices and hyperedges are distinct.
The \define{length} of a path $\gamma$ is the number $\mathrm{length}(\gamma)$ of hyperedges it contains.   
The \define{shortest path distance}   is the function $d_\SP\colon V\times V\to \mathbb{N}$ defined by
\[\dist_\SP(u,v)=\min_\gamma \length(\gamma) \, ,
\]
where the minimum is taken over all paths $\gamma$ between $u$ and $v$.

\end{definition}
\begin{definition}
Let $n$ be a positive integer. The \define{$n$-neighbourhood graph} of a hypergraph $\HG$ is the graph $H_n$ with vertex set given by that of $\HG$ and set of edges given by hyperedge paths of length $n$:
\[
H_n=\{\{u,v\}\mid \text{there exist a path of hyperedges of length $n$ between $u$ and $v$} \}\, .
\]
We call  $H_1$  the \define{graph associated} to the hypergraph $\HG$. 

\end{definition}

To extend Ollivier Ricci curvature to hypergraphs, various approaches can be considered. These approaches involve generalising probability measures and distances between them, where in some  extensions the cardinality of hyperedges are taken into account.  Similarly as for graphs, curvature can be defined for two arbitrary connected vertices $x, y$ or more generally for any sets of vertices, connected by one or more hyperedges.
Similarly as for graphs, we only consider $k(\mu_x, \mu_y)$ for vertices $x,y$ belonging to the same connected component of a hypergraphs. 
To generalise random walks to hypergraphs, we consider two approaches from  \cite{hypergraphcurvature}:

\begin{definition}
We define the following random walks on the set of vertices $V$ of an undirected hypergraph $\HG=(V,H)$:
\begin{itemize}
    \item \define{equal-nodes} (EN) random walk where we pick a neighbour $z$ of vertex $x$ uniformly at random: $\mu_x(z):= \frac{1}{|N(x)|}$ where $N(x)=\left\{ z : d_\SP(x,z)=1 \right\}$
    \item \define{equal-edges} (EE) random walk which first picks a hyperedge e with $|e| \geq 2$, then picks a vertex $z \in e-\{x\}$, both uniformly at random: $\mu_x(z):= \frac{1}{deg(x)-|\{x \in e:|e|=1\}|}\sum_{x,z \in e} \frac{1}{|e|-1}$
  
\end{itemize}

\end{definition}

After selecting one of these random walks we use the same definition of curvature as before: 

\begin{definition}
Let $x,y$ be two vertices in the same connected component of an undirected hypergraph $\HG=(V,H)$. Let $\{\mu_x\}_{x\in V}$ be one of the two random walks defined above. We set 
 \[\kappa(x,y):=1-\frac{W_1(\mu_x,\mu_y)}{d_\SP(x,y)}\, .\]
 When the probability measures $\{\mu_x\}_{x\in V}$ are defined via the EN random walk, we  call $\kappa(x,y)$ the   \define{EN-ORC} of $(x,y)$ and and when they are defined using the EE random walk, we call $\kappa(x,y)$ the \define{EE-ORC} of $(x,y)$. 

\end{definition}
Note that  equal-edges  random walks incorporate the cardinalities of hyperedges, whereas the equal-nodes random walk only records adjacency. \\

Throughout this section, we assume that the hypergraph $\HG=(V,H)$ has no loops,
i.e., all hyperedges satisfy $|e|\ge2$.

First,  we observe that weak equivalences can be derived in complete analogy with the graph case, allowing graph-based results to extend naturally to hypergraphs for both structural and regular equivalences. Note that $H_1$ captures all adjacency relations among the vertices of $\HG$. Therefore, whenever only the existence of connections matters, rather than their type, number, or cardinality, the equivalence relations in $\HG$ and in $H_1$ yield the same classes for both structural and regular notions. By the definition of weak equivalences, we have the following:

\begin{lemma}\label{L:H_1}
    Two vertices in $\HG$ are weakly structurally equivalent (resp.\ weakly regularly equivalent) iff they are structurally equivalent (resp. regularly equivalent) in $H_1$.
\end{lemma}

\begin{theorem}\label{T:hyper str}
    Let $\HG=(V,H)$ be an undirected hypergraph. Then: 
     \begin{enumerate}[label=(\roman*{})]
        \item Two vertices $x$ and $y$ are weakly structurally equivalent iff the EN-ORC of   $(x,y)$ is equal to $1$.
        \item If two vertices $x$ and $y$  are strongly structurally equivalent then the EE-ORC of   $(x,y)$ is equal to $1$. 
    \end{enumerate}
\end{theorem}
\begin{proof}
 To see the first case we recall that in weak equivalence $z$ is a neighbour of $x$  if and only if it is a neighbour of $y$ independent of the cardinality of the connecting hyperedge (hyperedges); we should show that EN-ORC of the pair $(x,y)$ is one. For EN-walks, for every such $z$, by definition $\mu_x(z)=\mu_y(z)$  and therefore $W_1(\mu_x,\mu_y)=0$ or equivalently $\kappa (x,y)=1$. For the reverse, if $\kappa (x,y)=1$, this means that $\mu_x=\mu_y$ and in particular they have exactly the same neighbours.
\\
For the second case,
if two vertices $x$ and $y$ are strongly structurally equivalent, then they are connected by hyperedges of the same cardinality to the same neighbours, and therefore their EE random walks are identical. This means that $W_1(\mu_x,\mu_y)=1$ and in particular their EE-ORC is equal to $1$.
\end{proof}

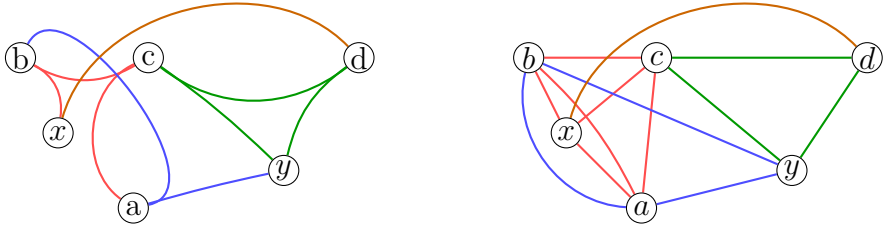
\begin{figure}[h]
\centering
\centering
\begin{tikzpicture}[scale=1, transform shape,
  vertex/.style={circle, draw=black, fill=white, minimum size=4mm, inner sep=0pt},
  edgeA/.style={thick, red!70, bend left=20},
  edgeB/.style={thick, blue!70, bend left=20},
  edgeC/.style={thick, green!60!black, bend left=20},
  edgeD/.style={thick, orange!80!black, bend left=20}
]
\node[vertex] (v1) at (-1,0) {$x$};
\node[vertex] (v2) at (0,-1) {a};
\node[vertex] (v3) at (2,-0.5) {$y$};
\node[vertex] (v4) at (-1.5,1) {b};
\node[vertex] (v5) at (0.2,1) {c};
\node[vertex] (v6) at (3,1) {d};

\begin{pgfonlayer}{background}

  \draw[edgeA] (v4) to[bend right=30] (v5);
  \draw[edgeA] (v2) to[bend left=60] (v5);
  \draw[edgeA] (v1) to[bend right=35] (v4);
  
  \draw[edgeB] (v2) to[bend left=2] (v3);
  \draw[edgeB] (v2) to[bend right=115] (v4);

  \draw[edgeC] (v5) to[bend right=38] (v6);
  \draw[edgeC] (v3) to[bend right=4] (v5);
  \draw[edgeC] (v6) to[bend right=20] (v3);
  \draw[edgeD] (v1) to[bend right=300] (v6);
\end{pgfonlayer}
\end{tikzpicture}
\qquad\qquad
\begin{tikzpicture}[scale=1, transform shape,
  vertex/.style={circle, draw=black, fill=white, minimum size=4mm, inner sep=0pt},
  edgeA/.style={thick, red!70, bend left=20},
  edgeB/.style={thick, blue!70, bend left=20},
  edgeC/.style={thick, green!60!black, bend left=20},
  edgeD/.style={thick, orange!80!black, bend left=20}
]

\node[vertex] (v1) at (-1,0) {$x$};
\node[vertex] (v2) at (0,-1) {$a$};
\node[vertex] (v3) at (2,-0.5) {$y$};
\node[vertex] (v4) at (-1.5,1) {$b$};
\node[vertex] (v5) at (0.2,1) {$c$};
\node[vertex] (v6) at (3,1) {$d$};

\begin{pgfonlayer}{background}

  \draw[edgeA] (v4) to[bend right=0] (v5);
  \draw[edgeA] (v4) to[bend right=0] (v1);
  \draw[edgeA] (v4) to[bend left=10] (v2);
  \draw[edgeA] (v5) to[bend right=0] (v2);
  \draw[edgeA] (v1) to[bend right=0] (v2);
  \draw[edgeA] (v1) to[bend right=0] (v5);
    
  \draw[edgeB] (v2) to[bend left=0] (v3);
  \draw[edgeB] (v2) to[bend left=50] (v4);
  \draw[edgeB] (v4) to[bend left=0] (v3);
  \draw[edgeC] (v5) to[bend left=0] (v6);
  \draw[edgeC] (v3) to[bend right=0] (v6);
  \draw[edgeC] (v5) to[bend right=0] (v3);
  \draw[edgeD] (v1) to[bend right=300] (v6);

\end{pgfonlayer}

\end{tikzpicture}
\caption{An undirected hypergraph (left) with four hyperedges colored by red, green, blue and orange, and its associated graph $H_1$ (right).}\label{F: ex 1}
\end{figure}

\begin{example} Consider the hypergraph in Figure \ref{F: ex 1} (left). The two vertices 
 $x$ and $y$  
 share the same neighbours, have the same degree and are weakly structurally equivalent, but fail to be strongly structurally equivalent. This can be seen by computing the curvature of $(x,y)$: we have that the EE-ORC of $(x,y)$ is different from $1$ and therefore by Theorem \ref{T:hyper str} (ii) $x$ and $y$ are not strongly structurally equivalent. On the other hand, the EN-ORC of $(x,y)$ is equal to $1$, and thus $x$ and $y$ are weakly structurally equivalent. For the latter, one can equivalently also consider the ORC of $(x,y)$ in the graph associated to $\HG$ (see Figure \ref{F: ex 1} (right)).
   \end{example}

\begin{figure}
\[
\begin{tikzpicture}[
  vertex/.style={circle, draw, minimum size=6mm},
  >={Stealth}
]

\node[vertex] (x) at (0.5,0) {$x$};
\node[vertex] (y) at (5,0) {$y$};
\node[vertex] (a) at (3,1) {$a$};
\node[vertex] (b) at (3,0) {$b$};
\node[vertex] (c) at (3,-1) {$c$};

\coordinate (split) at (1.5,0);

\draw[-, ] (x) -- (split);
\draw[- ] (split) -- (a);
\draw[-] (split) -- (b);
\draw[-] (split) -- (c);
\draw[-] (y) -- (a);
\draw[-] (y) -- (b);
\draw[-] (y) -- (c);
\end{tikzpicture}
\]
\caption{The two vertices $x$ and $y$ have the same neighbours, different degrees, and are weakly structurally equivalent. Moreover, we have that the EE-ORC of $(x,y)$ is equal to $1$, but $x$ and $y$ are not strongly structurally equivalent. }\label{F:counter geom char}
\end{figure}
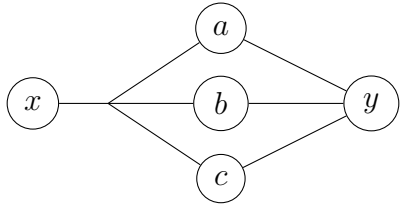

We note that the converse implication in item (ii) of Theorem \ref{T:hyper str} does not hold, and that   it is not possible to obtain a full characterisation of strong structural equivalences for hypergraphs through the two types of Ollivier Ricci curvature we consider here. This is easily seen by considering the example of two vertices that have the same three neighbours but are connected to them through a single hyperedge vs.\ three edges (see Fig. \ref{F:counter geom char}). In this example the three vertices $a,b,c$ have probability $1/3$ under both EE and EN random walks, for both $x$ and $y$.
On  one hand this raises the question of whether it could be possible  to consider a different type of random walk that allows to tell the vertices $x$ and $y$ apart. We discuss this further in Section \ref{S:conclusion}. On the other hand, we note that while we do not have a full characterisation, EE-ORC can still be a useful measure of approximate strong structural equivalence, as we discuss in Section \ref{S:approx str eq}.

Finally, we note that by imposing  further assumptions, we can obtain the following converse implication for item (ii) in Theorem \ref{T:hyper str}: 
 \begin{lemma}
Let $x$ and $y$ be two vertices of a hypergraph. 
Assume that
\begin{enumerate}
    \item for every common neighbour $z$ of $x$ and $y$, there exists exactly one hyperedge containing $\{x,z\}$ and exactly one hyperedge containing $\{y,z\}$;
    \item $d_x=d_y$.
\end{enumerate}
If\/ $\mathrm{EE\text{-}ORC}(x,y)=1$, then $x$ and $y$ are strongly structurally equivalent.
\end{lemma}
\begin{proof}
Since $\mathrm{EE\text{-}ORC}(x,y)=1$, we have $W(\mu_x,\mu_y)=0$, and hence
\[
\mu_x=\mu_y .
\]
In particular, $\mu_x$ and $\mu_y$ have the same support, which implies that they have the same neighbours. 
Let $z$ be an arbitrary common neighbour of $x$ and $y$. 
By assumption, there exists exactly one hyperedge $e_{xz}$ containing $\{x,z\}$ 
and exactly one hyperedge $e'_{yz}$ containing $\{y,z\}$. 
By the definition of the EE walk, we have
\[
\mu_x(z)=\frac{1}{d_x}\cdot\frac{1}{|e_{xz}|-1},
\qquad
\mu_y(z)=\frac{1}{d_y}\cdot\frac{1}{|e'_{yz}|-1}.
\]
Since $\mu_x=\mu_y$ and $d_x=d_y$, it follows that
\[
|e_{xz}| = |e'_{yz}|.
\]
Because $z$ was arbitrary, this equality holds for every common neighbour $z$ of $x$ and $y$. 
Therefore, $x$ and $y$ are connected to the same neighbours via hyperedges of the same cardinality, 
which means that $x$ and $y$ are strongly structurally equivalent.
\end{proof}

\begin{remark}[Hierarchy of structural equivalence notions in hypergraphs]\label{R:hierarchy}
The proof of the above lemma relies on rather strong local assumptions,
most notably the uniqueness of hyperedges connecting a vertex to each of its neighbours and equality of degrees. 

These assumptions are sufficient to ensure that equality of the EE measures
$\mu_x=\mu_y$ forces a complete matching of local hyperedge cardinalities, but they are
not intrinsic to the definition of the EE random walk or to EE-ORC.
This motivates the introduction of a hierarchy of structural equivalence notions of
increasing strength, reflecting different levels of combinatorial information that can be
recovered from $\mu_x$ and from curvature-based quantities. We note  in all  the following notions of equivalence we have that  the neighbours of $x$ and $y$ are identical, which means that all these notions  are stronger than weak structural equivalence.
\medskip

\noindent\textbf{(i) EE-measure equivalence.}
Two vertices $x$ and $y$ are said to be \emph{EE-measure equivalent} if
\[
\mu_x=\mu_y .
\]
This is the weakest notion of equivalence considered here and is exactly characterized by
the condition $\mathrm{EE\text{-}ORC}(x,y)=1$.
It captures equivalence with respect to one-step EE diffusion, but does not encode the
underlying hyperedge incidence structure uniquely.

\medskip
\noindent\textbf{(ii) Weighted-neighbourhood equivalence.}
Vertices $x$ and $y$ are said to be \emph{weighted-neighbourhood equivalent} if, for every
vertex $z$,
\[
\sum_{\substack{x,z\in e}}\frac{1}{|e|-1}
=
\sum_{\substack{y,z\in e}}\frac{1}{|e|-1}.
\]
We refer to this notion as weighted-neighbourhood equivalence because it compares vertices
by the total hyperedge-induced weights assigned to each neighbour, obtained by aggregating
contributions of the form $(|e|-1)^{-1}$ from all incident hyperedges.
Under the loop-less assumption, when $d_x=d_y$, this notion coincides with EE-measure
equivalence by definition of the EE random walk.

\medskip
\noindent\textbf{(iii) Multiset structural equivalence.}
Vertices $x$ and $y$ are \emph{multiset structurally equivalent} if
they have the same neighbours and, for every common neighbour $z$, the multisets
\[
\bigl\{|e|:\ e\in E,\ \{x,z\}\subseteq e\bigr\}
\quad\text{and}\quad
\bigl\{|e|:\ e\in E,\ \{y,z\}\subseteq e\bigr\}
\]
coincide.
This notion preserves the full spectrum of hyperedge cardinalities incident to each
neighbour, while allowing multiple hyperedges and without imposing a bijection between
individual hyperedges.
It is strictly stronger than EE-measure equivalence, but weaker than full combinatorial
incidence matching.

\medskip
\noindent\textbf{(iv) Strong  structural equivalence.}
We note that in this notion that we originally defined, the 
local hyperedge incidence patterns  of $x$ and $y$ agree in a fully combinatorial sense, in particular
matching hyperedges (with multiplicity) connecting $x$ and $y$ to each common neighbour and
preserving cardinalities.
This is the strongest notion considered here and corresponds to a local isomorphism of
incidence stars.

\medskip
\noindent\textbf{Relations between the notions.}
The above definitions form a natural implication chain:
\[
\text{strong eq}
\;\Longrightarrow\;
\text{multiset eq}
\;\Longrightarrow\;
\text{weighted-neighbourhood eq}
\;\Longrightarrow\;
\text{EE-measure eq}.
\]
In general, none of the reverse implications holds.
The above lemma shows that, under additional assumptions such as unique
hyperedge incidence and same degrees, EE-measure equivalence is
sufficient to recover strong  structural equivalence.
Without these assumptions, curvature-based quantities cannot distinguish between
different local incidence patterns that induce the same EE transition measures.
\\
We do not investigate further properties of these notions here and leave a systematic study for future work.
\end{remark}

We conclude this section by giving two direct generalisations of results for graphs to hypergraphs. First, we have the following results based on considering $H_2$:

\begin{corollary}\label{C:H_2 cc}
   Let $\HG$  be an undirected hypergraph. Then we have:
  \begin{enumerate}[label=(\roman*{})]
    \item  Weak structural classes of $\HG$ are vertex sets of complete subgraphs of  connected components of $H_2$.
   
    \item  Complete subgraphs of connected components of $H_2$ induce weak structural classes of $\HG$ if and only if the EN-ORC curvature of every pair of vertices  in the same subgraph is equal to $1$ in $\HG$, and the subgraph is the maximal subgraph with respect to this property. 

    \item \label{I: cc H_2} Connected components of $H_2$ induce weak regular classes of $\HG$. Moreover, if $H_2$ has two connected components, these components induce a weak regular partition that is not properly contained in any non-trivial weak regular partition.  
    \item \label{I: scc H_2}  Complete subgraphs of the connected components of $H_2$ such that the EN-ORC  of all  the pairs of the vertices in each subgraph is equal to $1$ in $\HG$ induce weak regular classes of $\HG$.
    \end{enumerate}
    
\end{corollary}
\begin{proof}
The proof is a direct consequence of   Lemma \ref{L:H_1} and Theorem \ref{geom char eq} (for (i) and (ii)), Theorem \ref{T: char cc G_2} (for iii), Corollary \ref{T:subc} (for (iv)).
\end{proof}
 
The following is a direct generalisation of Theorem \ref{T: reg G_n} to hypergraphs:

\begin{corollary}\label{C: H_n reg}
    Let $\HG$  be an undirected hypergraph.   If each vertex has at least two neighbours, then the connected components of $H_n$ induce weak regular classes of $\HG$.
\end{corollary}
\begin{proof}
    The proof is a direct consequence of Lemma  Lemma \ref{L:H_1} and Theorem \ref{T: reg G_n}.
\end{proof}

\section{Approximate structural equivalences for graphs and hypergraphs}\label{S:approx str eq}

As we saw, structural equivalence   is a very restrictive notion. Thus, in practice, one  works with measures of approximate structural equivalence, such as cosine similarity, see Section \ref{D:cosine sim}.
A different measure of approximate structural equivalence is suggested by our geometric characterisation, not only for graphs but also for hypergraphs: we know that two vertices $x,y$ of a graph are structurally equivalent if and only if the Ollivier-Ricci curvature of $(x,y)$ is equal to $1$. This result extends directly to hypergraphs: weak structural equivalence can be characterised via EN-ORC, and EE-measure equivalence via EE-ORC (see Theorem~\ref{T:hyper str} and Remark~\ref{R:hierarchy}). By relaxing these equalities, we obtain natural measures for approximating structural equivalence: the closer the corresponding curvature is to $1$, the more structurally equivalent the two vertices are. In this section, we aim to better understand the type of information that Ollivier--Ricci curvature provides for both undirected graphs and hypergraphs.

We first make some observations.

\begin{itemize}
    \item The curvature $\kappa(x,y)$ is equal to $1$ if and only if  $x,y$ are structurally equivalent. We note that in this case $d_\SP(x,y)\neq1$.  Similarly, for vertices in hypergraphs, this holds with EN-ORC for weak structural equivalence and with EE-ORC for EE-measure equivalence (see Remark \ref{R:hierarchy} for the latter).  
    \item In both graphs and hypergraphs, if $d_\SP(x,y)>2$, then $x$ and $y$ have no neighbours in common and we call them \define{fully dissimilar}.  
    \item In a graph, if  $ d_\SP(x,y)=2$ : $-1 < \kappa(x,y)\leq1 $. Therefore the closer the curvature of this pair is to one, the more structurally similar they are and the closer the curvature is to $-1$, they are less and less similar. The same statement holds for vertices $x$ and $y$ in a hypergraph, with weak (respectively EE-measure) equivalence based on EN-ORC (respectively EE-ORC).
    \item If $ d_\SP(x,y)=1$ : $-2 < \kappa(x,y) < 1 $; similar to the previous case, the closer the curvature is to one (in this case it can never be one), the more similar the two vertices are but in opposite the closer this number is to $-2$ , they are less and less similar. Replacing $\kappa(x,y)$ by EN-ORC or EE-ORC yields analogous conclusions for weak and EE-measure equivalences in hypergraphs, respectively.
\end{itemize}

Thus, for vertices at distance at most $2$, the curvature is a number between $-2$ and $1$. On the other hand, cosine similarity is a number between $0$ and $1$. We note that  cosine similarity has been introduced for graphs, and a direct generalisation to hypergraphs would in general provide a different type of measure, since for hypergraphs the degree of a vertex is in general smaller than the number of its neighbouring vertices. For hypergraphs, we  believe that what is of most interest is considering directly measures of approximate equivalences given by the random walk equivalence notions discussed in Remark \ref{R:hierarchy}. We leave the exploration of these random walk equivalence notions for future work.  

How are curvature and cosine similarity related? 
To obtain the connection between the two notions in graphs, let us suppose we want to move sand $\mu_x$ to holes $\mu_y$ in an optimal way. 
We note that, depending on the distance between $x$ and $y$, some sand and holes  may already coincide;  as mentioned in \cite{jost2013, eidi2020}, this amount is  equal to $\frac{\eta_{xy}}{d_x \vee d_y}$ where $d_x \vee d_y$ is the maximum of the degrees of $x$ and $y$.
Importantly, there always exists an optimal transportation plan that does not move  the sand already in holes \footnote{For graphs this is a  consequence of the triangle inequality, see  \cite[Lemma 4.1]{BCLMP18}. More generally, as noted in \cite{BCLMP18}, this follows from the invariance of Kantorovich-Rubinstein distance under mass subtraction \cite[Corollary 1.16]{villani03}.}. 
As observed  in \cite{eidi2020}, for two neighbouring vertices $x$ and $y$ one has $W_1(\mu_x,\mu_y)= \sum_{i=0}^{3} i \mu_i  $ where $\mu_i$ is the amount of sand that should be moved with distance $i$, in an optimal transportation plan; if $x$ and $y$ are not neighbours we can simply have the same formula by extending the range of this summation to $d(x,y)+2$ as that is the maximum distance of sand from a hole.
Since there exists an optimal transportation plan that leaves the sand already in holes untouched, for such a transportation plan we have $\mu_0=\frac{\eta_{xy}}{d_x \vee d_y}$.
This  observation allows us to obtain the following relationship between cosine similarity and Ollivier Ricci curvature:

\begin{proposition}\label{P:bounds curvature}
Let $x,y$ be any pair of distinct vertices in a graph. We  have:
\[
 3\frac{\sigma_{xy}\sqrt{d_xd_y}}{d_x\vee d_y} -2 \leq \kappa(x,y) \leq  1-\frac{1-\sigma_{xy}}{d(x,y)}\, .
\]
If we instantiate these inequalities for different values of cosine similarity, we obtain $\kappa(x,y)=1$ for $\sigma_{xy}=1$, and furthermore:
\begin{itemize}
\item if $\sigma_{xy}=0$ then we have $-2\leq \kappa (x,y)\leq 1- \frac{1}{d(x,y)}$,
\item if $\sigma_{xy}=s$ with $0<s<1$ then  we have $ 3s\frac{\sqrt{d_xd_y}}{dx\vee d_y} -2\leq \kappa (x,y)\leq 1- \frac{1-s}{d(x,y)}\, .$
\end{itemize}
\end{proposition}

\begin{proof}
We first prove the upper bound.  Let  $\mu_i$ denote the amount of mass that is moved across distance $i$  in an optimal transportation plan. We  thus have 
\[
W_1(\mu_x,\mu_y)=\sum_{i>0}i \, \mu_i  \geq \sum_{i>0}\mu_i = 1-\mu_0\, ,
\]
since $\sum_{i\geq 0}\mu_i=1$. 
Without loss of generality, we can   assume that $\mu_0$ equals  the amount of mass in common between the two probability measures, since there  exists a transportation plan leaving this quantity untouched.

Thus, 
\begin{alignat}{2}W_1(x,y)&\notag \geq 1-\mu_0\\
&\notag =1-\frac{\eta_{xy}}{d_x\vee d_y}\\
&\notag =1-\frac{\sigma_{xy}\sqrt{d_xd_y}}{d_x\vee d_y}\\
&\notag \geq 1-\sigma_{xy}\, ,
\end{alignat}
where the last inequality holds because $\frac{\sqrt{d_xd_y}}{d_x\vee d_y}\leq 1$. Thus, we obtain
\[
\kappa(x,y)=1-\frac{W_1(\mu_x,\mu_y)}{d(x,y)} \leq 1-\frac{1-\sigma_{xy}}{d(x,y)} \, .
\]\\

For the lower bound, 
we first observe: 
\[
W_1(\mu_x,\mu_y)=\sum_{i=1}^{d(x,y)+2}i \,\mu_i \leq (d(x,y)+2)\sum_{i=1}^{d(x,y)+2}\mu_i=(d(x,y)+2)(1-\mu_0)\, . \]

We thus get the following lower bound:
\begin{alignat}{2}
\kappa(x,y)=1-\frac{W_1(\mu_x,\mu_y)}{d(x,y)}&\notag \geq 1 - \frac{(d(x,y)+2)(1-\mu_0))}{d(x,y)}\\
&\notag = 1 - (d(x,y)+2)\frac{(1-\frac{\sigma_{xy}\sqrt{d_xd_y}}{dx\vee d_y})}{d(x,y)} \\
&\notag = \frac{\frac{\sigma_{xy}\sqrt{d_xd_y}}{dx\vee d_y}(d(x,y)+2)-2}{d(x,y)}\, .
\end{alignat}
Furthermore, we remark that if $\kappa(x,y)\geq \kappa_0$ for every edge $\{x,y\}$ then $\kappa(x,y)\geq \kappa_0$ for any pair of vertices (this is a direct consequence of the triangle inequality for $W_1$).
We can thus take $d(x,y)=1$ and we obtain the lower bound
\[
\kappa(x,y)\geq 3\frac{\sigma_{xy}\sqrt{d_xd_y}}{dx\vee d_y} -2
\]
for any $x\ne y\in V$.

Finally, we  note that if $\sigma_{xy}=1$ then  we must have $d_x=d_y$ and $d(x,y)=2$. Thus, the lower bound simplifies to $1$ and same for the upper bound, and we obtain $1\leq \kappa(x,y)\leq 1$.
\end{proof}

A further  way to interpret the statement of  Proposition \ref{P:bounds curvature} is that it gives easily-derivable bounds on the curvature of two vertices, which  can be understood in very simple terms, since they are simply given by the proportion of common neighbours. For instance, if  two vertices are fully dissimilar, then we can say that the curvature is at most $1-\frac{1}{d(x,y)}$. In particular, for neighbour vertices we must have that the curvature is zero or negative.

\section{Conclusion and open questions}\label{S:conclusion}

In this work we have linked Ollivier Ricci curvature and random walks on graphs and hypergraphs with the notions of structural and regular equivalence. We believe that our work is only a first step in this direction, and that the connection established raises many interesting questions. One of these concerns random walks on hypergraphs. Our study of characterisations of structural equivalences for hypergraphs through curvature highlighted the fact that current existing random walk notions on hypergraphs are not able to discern local adjacency patterns in hypergraphs; for instance, all existing notions of random walks would assign the same values to the nodes  $x$ and $y$ in Figure \ref{F:counter geom char}. We can envision applications in which this distinction might be of importance, for instance in  social systems modelling, where a hyperedge might represent coattendance of events, and the random walk could encode the probability that a person might talk with another person during the same event. In such a scenario, the situations of the two nodes $x$ and $y$ are drastically different, and thus should be better captured by a different type of random walk.



In extending the curvature characterisation to hypergraphs, we have noted that different types of equivalences based on random walks are better suited at capturing different local incidence patterns of similarity. We reserve the systematic study of these novel notions of equivalences for hypergraphs to future work. 


Finally, we note that regular equivalences and  Ollivier-Ricci curvature    have recently  found application in machine learning;  regular equivalences have been used in representation learning \cite{TCWYZ18}, while  Ollivier-Ricci curvature has found application in the study of bottlenecks and oversquashing in graph neural networks \cite{2025arXiv250606582T}. This thus raises the question of what new insights the connection given in this work may give for problems in  geometric deep learning. 

\section{Acknowledgements}
We acknowledge funding from the SALTO exchange programme between the Max Planck Gesellschaft (MPG) and the Centre National de la Recherche Scientifique (CNRS), and the Fondation Mathématiques Jacques Hadamard.

\bibliography{main}
\bibliographystyle{alpha}
\end{document}